\documentclass[11pt]{article}
\usepackage{amsfonts}
\usepackage{amssymb}
\usepackage{amsthm}
\usepackage{amsmath}
\usepackage{graphicx}
\usepackage{indentfirst}
\usepackage{cite}
\usepackage{mathrsfs}
\usepackage{graphics}

\usepackage{amsfonts,amsmath,amstext,amsopn,amsbsy,amscd}
\usepackage{amsxtra,amssymb,latexsym,mathrsfs}
\usepackage{eufrak}
\usepackage{epsfig}
\usepackage{pstricks,multido,pst-node}

\newtheorem{theorem}{\textbf{Theorem}}[section]
\newtheorem{lemma}{\textbf{Lemma}}[section]
\newtheorem{proposition}{\textbf{Proposition}}[section]
\newtheorem{corollary}{\textbf{Corollary}}[section]
\newtheorem{remark}{\textbf{Remark}}[section]
\newtheorem{definition}{\textbf{Definition}}[section]

\allowdisplaybreaks[4]

\def\be{\begin{equation}}
\def\ee{\end{equation}}
\def\bea{\begin{eqnarray}}
\def\eea{\end{eqnarray}}
\def\bt{\begin{theorem}}
\def\et{\end{theorem}}
\def\bl{\begin{lemma}}
\def\el{\end{lemma}}
\def\br{\begin{remark}}
\def\er{\end{remark}}
\def\bp{\begin{proposition}}
\def\ep{\end{proposition}}
\def\bc{\begin{corollary}}
\def\ec{\end{corollary}}
\def\bd{\begin{definition}}
\def\ed{\end{definition}}

\def\non{\nonumber }

\begin{document}

\title{Diffusion limit of kinetic equations for multiple species charged particles}

\author{
{\sc Hao Wu} \footnote{School of Mathematical Sciences and Shanghai Key Laboratory for Contemporary Applied Mathematics, Fudan University,
200433 Shanghai, China, Email: \textit{haowufd@yahoo.com}},\ \
 {\sc Tai-Chia Lin} \footnote{Department of Mathematics, National Taiwan University, No.1, Sec. 4,
Roosevelt Road, Taipei 106, Taiwan, Email:
\textit{tclin@math.ntu.edu.tw}},\ \ and
{\sc Chun Liu}
\footnote{Department of Mathematics, Penn State University, State
College, PA 16802, Email: \textit{liu@math.psu.edu}}
}

\date{\today}

\maketitle


\begin{abstract}
In ionic solutions, there are multi-species charged particles (ions) with different properties like mass, charge etc. Macroscopic continuum models like the Poisson--Nernst--Planck (PNP) systems have been extensively used to describe the transport and distribution of ionic species in the solvent. Starting from the kinetic theory for the ion transport, we study a Vlasov--Poisson--Fokker--Planck (VPFP) system in a bounded domain with reflection boundary conditions for charge distributions and prove that the global renormalized solutions of the VPFP system converge to the global weak solutions of the PNP system, as the small parameter related to the scaled thermal velocity and mean free path tends to zero. Our results may justify the PNP system as a macroscopic model for the transport of multi-species ions in dilute solutions.
\\

\end{abstract}

\section{Introduction}
\setcounter{equation}{0}

The transport of ions in different biological environments is very important in our life and it has attracted more and more attentions recently \cite{Hille, HH, Eisen,NE}. In biological problems, the ionic solutions usually consist of charged particles (ions) like sodium Na$^+$, potassium K$^+$, calcium Ca$^{2+}$ and chloride Cl$^-$ etc, which have different but comparable masses, charge valencies and sizes. These differences have dramatic effects on the dynamics of multi-species ions which produce the functions of cells in biological system, e.g., the ion channels. To study the dynamics of multi-species ions, molecular dynamics simulations (MD) using microscopic models (from Newton's laws) to describe charge particle trajectories are popular and useful but expensive because the models are usually huge and the numerical computation time is very long (see, for instance,~\cite[Chapter~10]{Hille}).

To increase the efficiency of numerical simulations, one may use the (multi-species) Poisson--Nernst--Planck (PNP) system~\cite{Nerst, Co, Eisen3, Eisen4, Liu14}, which is a macroscopic model to describe multi-species ion transport. Conventionally, the PNP system consists of coupled diffusion-convection equations and the Poisson equation being represented as follows:
\be
\begin{cases}
 \partial_t c_i=\nabla \cdot J_i,\\
 J_i=d_i \left(\nabla c_i+ \displaystyle{\frac{q_i}{k_B T}} c_i \nabla \phi\right),\\
 -\nabla \cdot( \epsilon\nabla \phi)=\sum_{i=1}^N q_i c_i+D(x).
 \end{cases} \label{pnpp}
\ee
Here, $c_i$ ($i=1,2,...,N$) stand for densities of
charged particles in the ionic solution and $\phi$ is the self-consistent electric potential.  Besides, $q_i$ are the (positive or negative) charges of particles, $J_i$ are the ionic flux densities, $d_i$ are their diffusion coefficients, $\epsilon$ is the dielectric
coefficient, $k_B$ is the Boltzmann constant, $T$ is the temperature and $D(x)$ is the permanent charge density in the domain. The PNP system \eqref{pnpp} is also one of the fundamental macroscopic models in the study of transport of carriers in semiconductors, see, e.g., \cite{Ma,MRS,J,Ga}. Concerning the mathematical analysis, the initial value problem and the initial boundary value problem of the PNP system have been extensively studied in the literature, we refer to \cite{AMT,BD,BH,FI,Ga,GaG, MRS,JA} and the reference cited therein.

The PNP system \eqref{pnpp} provides a continuum description of the evolution of charged particles via macroscopic (averaged) quantities, e.g., the particle density, the current density etc., which have cheaper costs for numerics. Such continuum models can be (formally) derived from kinetic models by coarse graining methods, like the moment method, the Hilbert expansion method and so on \cite{JA,MRS,M}. Although many results for the PNP systems have been obtained, it seems that none of them can reveal basic principles like gating and selectivity of ion channels. Recently, new PNP type systems have been derived and the selectivity of ion channels have been simulated successfully~\cite{Eisen3, HLLE12, HEL11, HFEL, LE13,XSL13}. In order to justify these continuum models, here we develop the kinetic theory for the PNP system like \eqref{pnpp} as the first step work. Our goal in the present paper is to rigourously justify the PNP system for dilute ionic solutions consisting of multi-species charged particles, by studying the diffusion limit of a suitable kinetic system. We will continue to study the kinetic theory for those new PNP type systems as in \cite{Eisen3, HLLE12, HEL11, HFEL, LE13,XSL13} in the near future.

In this paper, we consider the case that the motion of multi-species charged particles is governed by the electrostatic force coming from their (self-consistent) Coulomb interaction. We also assume that the momentum of charged particles with collision is small and ignorable. Then the collision term in the kinetic equation may be approximated by the Fokker--Planck operator that describes the Brownian force~\cite{CH}, and the resulting kinetic system becomes the Vlasov--Poisson--Fokker--Planck (VPFP) system as follows:
\bea
 &&\partial_t f_i +v\cdot\nabla_x f_i-\frac{z_i q}{m_i}\nabla_x\phi \cdot\nabla_v f_i=\frac{1}{\tau_i}\mathcal{L}^i_{FP}(f_i),\label{f1}\\
 && -\epsilon_0\Delta_x \phi= q\left(\sum_{i=1}^N z_i \int_{\mathbb{R}^d} f_i dv+D(x)\right),\label{P1}
 \eea
 where $\mathcal{L}^i_{FP}$ $(i=1,...,N)$ are the Fokker--Planck operators such that
 \bea
  \mathcal{L}^i_{FP}(f_i)=\nabla_v\cdot\left( vf_i+  \theta_i\nabla_v f_i\right).\non
 \eea
Here, $\epsilon_0>0$ is the vacuum permittivity, $q>0$ is the positive elementary charge. For $i=1,...,N$, the state of each species is given by a distribution function $f_i(t, x, v)\geq 0$, i.e., a probability density in the $(x, v)$-phase space at time $t$ ($f_i dxdv$ is the number
of the $i$-th species charged particles at time $t$ located at a volume element $dx$ about the position $x$ and having
velocities in a volume $dv$ about the value $v$). Besides, $z_i\in\mathbb{Z}$ are the valencies for the $N$-species charged particles, $m_i$ are the masses, $\tau_i$ are relaxation time due to collisions of the particles with the thermal bath, $\sqrt{\theta_i}$ are the thermal velocities given by $\sqrt{\theta_i}=\sqrt{2k_B T_b m_i^{-1}}$ and  $T_b$ is the temperature of the thermal bath.

In plasma physics, the VPFP system \eqref{f1}--\eqref{P1} with $N=1$ (i.e., the single species case) is reasonable because the mass ratio between the ions and electrons is huge, only the evolution of the electrons is described in terms of a distribution function in the resulting system, and those `heavy' ions are supposed to be static. For such a case, the existence and uniqueness of solutions to the initial value problem or the initial boundary value problem of the VPFP system have been investigated in the literature. We refer to \cite{B93,VO,RW} for results on the classical solutions and to \cite{B95,GaS,Ca98,V} for weak solutions and their regularity. Concerning the long-time behavior of the VPFP system, we refer to \cite{BCS97,BD95,CSV}. Instead of the single species case with $N=1$, here we study the system \eqref{f1}--\eqref{P1} with $N\geq 2$ for multi-species charged particles, which is more complicated due to the (nonlocal) interactions between particles via the Poisson equation \eqref{P1} (i.e., the Coulomb interaction).

Suitable scalings of the VPFP system should be introduced in order to study its diffusion limit.  Let $L$ be the characteristic length. We denote by $N_0$ the characteristic value for the concentration of particles and by $\Phi_0$ the characteristic variation of the electric potential over $L$. Since we have to treat mutiple species of charged particles that have different masses and charges, it is convenient to introduce a `reference particle' with mass $m_{ref}$, electric charge $z_{ref} q$ (with $z_{ref}=1$), relaxation time $\tau_{ref}$ and thermal velocity $\theta_{ref}$. The microscopic variation as well as the drift velocity for the reference particle are given by $V_{ref}=\sqrt{\theta_{ref}}$, $ U_{ref}=\tau_{ref} \frac{q}{m_{ref}}\frac{\Phi_0}{L}$, respectively.  Choosing  the following scaling (with respect to the reference particle) $t\to T_0 t'$, $ x\to L x'$, $ v\to V_{ref} v'$, $ T_0=\frac{L}{U_{ref}}$ and the change of unknowns
 $ f_i(t,x,v)= N_0 V_{ref}^{-d}f_i'(t',x',v')$, $\phi(t,x,v)=\Phi_0 \phi'(t',x',v')$, $D(x)=N_0 D'(x')$,
 we obtain the rescaled VPFP equations (drop the prime for simplicity):
 \bea
 && \partial_tf_{i}+\nu v\cdot \nabla_{x} f_i-\frac{\kappa_i z_i}{\varepsilon} \nabla_{x}\phi\cdot\nabla_{v} f_i = \frac{\zeta_i \nu}{\varepsilon}\nabla_{v}\cdot (vf_i+\kappa_i \nabla_{v}f_i), \quad i=1,...,N,\non\\
  && -\varpi \Delta_x \phi=\sum_{i=1}^N z_i\int_{\mathbb{R}^d} f_i dv+D(x),\non
 \eea
 where the dimensionless parameters $\nu$ (the `scaled' thermal velocity), $\varepsilon$ (the `scaled' thermal mean free path),  $\varpi$ and the ratios $\kappa_i$, $\zeta_i$ are given by
 \be
 \nu= \frac{V_{ref}}{U_{ref}}, \quad \varepsilon=\frac{\tau_{ref} V_{ref}}{L}, \quad \varpi =\frac{\epsilon_0 \Phi_0}{qN_0L^2}, \quad \kappa_i=\frac{m_{ref}}{m_i}, \quad \zeta_i=\frac{\tau_{ref}}{\tau_i}.\non
 \ee

 The case we are interested in this paper is called the low field limit (or the parabolic limit), which means that the drift velocity is small comparing with the thermal velocity, while the thermal velocity is small comparing to the relaxation velocity, and the two ratios have the same order of magnitude (cf. \cite{PS00, G05, ACGS01, GhMa10}):
  $$\nu\simeq\varepsilon^{-1} \quad \text{and}\quad \varepsilon<<1.$$
  For $\varepsilon>0$, taking $\nu =\varepsilon^{-1}$ (just for the sake of simplicity), we arrive at the rescaled VPFP system under low field scaling, which will be investigated in the remaining part of this paper:
    \bea
 &&\partial_t f_i^\varepsilon +\frac{1}{\varepsilon}v\cdot\nabla_x f_i^\varepsilon-\frac{\kappa_i z_i}{\varepsilon}\nabla_x\phi^\varepsilon\cdot\nabla_v f_i^\varepsilon
 =\frac{\zeta_i}{\varepsilon^2}L^i_{FP}(f_i^\varepsilon),\label{f}\\
 && -\varpi \Delta_x \phi^\varepsilon= \sum_{i=1}^N z_i \int_{\mathbb{R}^d} f_i^\varepsilon(t,x,v) dv+D(x),\label{P}
 \eea
 where the rescaled Fokker--Planck operators are given by
 \be
 L_{FP}^i (f_i^\varepsilon)= \nabla_{v}\cdot (v f_i^\varepsilon+ \kappa_i \nabla_{v}f_i^\varepsilon).\label{FKO}
 \ee

We recall that the diffusion limit of the VPFP system has been studied extensively in the literature (cf. \cite{GhMa10,G05, GNPS,NPS01,PS00} and the references therein).
In \cite{GhMa10, PS00, G05}, the authors studied the low field limit and proved the convergence of suitable solutions to the single species VPFP system towards a solution to the drift--diffusion--Poisson model in the whole space. In \cite{PS00}, under a suitable regularity assumption on the initial data, the convergence result was obtained globally in time in two dimensions and locally in time for the three dimensional case.  Later, the author proved in \cite{G05} a global convergence result in the two dimensional case, without any restriction on the time interval and the assumptions on the initial data were weakened with bounds only on the associated entropy and energy. Quite recently, in \cite{GhMa10} the authors established a global convergence result, without any restriction on the time interval or on the spatial dimensions, by working with the renormalized solutions (or free energy solutions, cf. \cite{DL88, Do99}). As pointed out in \cite{GhMa10}, the notion of renormalized solutions is natural for the problem, because the free energy of the VPFP system seems to be the only quantity that is uniformly bounded with respect to the small parameter $\varepsilon$ (i.e., the `scaled' mean free path). Even one works with more regular initial data such that the solutions can be defined in the usual weak sense without the need of renormalizing, one still has to use renormalization techniques to pass to the limit as $\varepsilon \to 0$. Besides, the use of
renormalization techniques together with an averaging lemma helps to remove the restriction on spatial dimensions and treat the nonlinear term $\nabla_x \phi \cdot \nabla_v f$, where the main difficulty comes from  (we refer to \cite{GhMa10} for more details).

 In this paper, we rigorously prove that for the multi-species case, the VPFP system~\eqref{f}--\eqref{P} converges to a rescaled PNP system as  $\varepsilon$ tends to zero in the low field limit. We generalized the techniques introduced in the previous works  \cite{GhMa10, Mi10, MT07}, to the case involving
multiple species of charged particles in a bounded region with reflection boundary conditions \cite{Mi10,Ca98,CGL01}. The specific boundary conditions recover the classical  no-flux boundary conditions of the PNP system. Different from the single species case in the literature, the previous arguments have to be modified in order to deal with the nonlocal interactions between different species of particles through the Poisson equation for the electric potential $\phi$. Besides, in order to deal with the integrals on the boundary, we shall make use of the Darroz\`es--Guiraud information \cite{DG66}, which helps to obtain the energy dissipation. Finally, effects of different but comparable quantities like masses and valencies of the charged particles will become obvious in our mathematical analysis.

Our results support the PNP system \eqref{pnpp} as a suitable model for  multi-species charged particles in dilute solution. As we mentioned before, several variants of the PNP system \eqref{pnpp} have recently been derived by using the energetic variational approaches \cite{HKL} to model important physical ingredients such as size (steric) effects for non-diluted solutions (cf. e.g., \cite{Eisen3,XSL13,LE13,HEL11,HLLE12,HFEL}) that are crucial in the study of the selectivity of ion channels in cell membranes \cite{Co, Ku, Hille}.  The total energy for these modified PNP systems consists of the entropic energy induced by the Brownian motion of ions, the electrostatic potential energy representing the coulomb interaction between the charged ions, and in particular, the repulsive potential energy caused by the excluded volume effect (e.g., the singular Lennard--Jones potential). Our result can be viewed as a starting point for the further investigation on the case of crowded ions. It would be interesting to study the diffusion-limit of suitable kinetic systems to obtain the modified PNP systems~\cite{Eisen3,HEL11}.

The remaining part of this paper is organized as follows. In Section 2, we  present the definition of renormalized solutions and state the main result on the diffusion limit of the VFFP system \eqref{f}--\eqref{P} (Theorem \ref{limit}). In Section 3, after deriving the energy dissipation of the VFFP system in bounded domain (Proposition \ref{energy}), which yields the necessary uniform estimates (Lemmas \ref{unie}, \ref{bdnphi}, \ref{r}), we proceed to prove our main result by using the renormalization techniques.

 \section{Preliminaries and main result}
 \setcounter{equation}{0}
 \subsection{Boundary and initial conditions}
Let $\Omega\subset \mathbb{R}^d$ $(d\geq 2)$ be a sufficiently smooth bounded domain. For instance, the outward
unit normal vector $\mathbf{n}(x)$ at $x \in \partial \Omega$ satisfies $\mathbf{n}\in W^{2,\infty}(\Omega, \mathbb{R}^d)$. The Lebesgue surface measure on $\partial \Omega$ will be denoted by $dS$.

Then we introduce the boundary conditions for the distribution functions. As in Cercignani's work \cite{Ce90, CLL95, CGL01} (see also \cite{BCS97, Mi10}), we define the sets of outgoing $(\Sigma^x_+)$ and incoming $(\Sigma^x_-)$ velocities at point $x \in\partial\Omega$ such that
$\Sigma^x_\pm := \{v\in \mathbb{R}^d: \ \pm v \cdot \mathbf{n}(x) > 0\}$ and denote the boundary sets
 $\Sigma_\pm =\{(x,v): \ x\in \Omega, \ v\in \Sigma^x_\pm\}$. Let $\gamma h$ be the trace of function $h$ and $\gamma_\pm h=\mathbf{1}_{(0,+\infty)\times \Sigma_\pm}\gamma h$.
Reflection boundary conditions for the kinetic equations take the form of integral (balance) relations between the densities of the particles on the outgoing and incoming velocity subsets of the boundary $\partial \Omega$ at a given time \cite{Ce90, CLL95, CGL01}. For instance, given $x\in \partial\Omega$ and $t>0$, we have (cf. \cite{BCS97}):
 \be
 \gamma_- f(t,x,v)=\int_{\Sigma^x_+}R(t,x; v,v^*) \gamma_+f(t,x,v^*) dv^*,\quad v\in \Sigma_-^x,\label{bcc}
 \ee
where $R$ represents the probability that a particle with velocity $v^*$ at time $t$  striking the boundary on $x$ reemerges at the same instant and location with velocity $v$.
If we consider $v'=-v$ for any $v\in \Sigma_-^x$ and take $R(t,x; v,v^*)=\delta_{v'}$ being the Dirac measure centered at $v^*=v'$, then we have $\gamma_-f(t,x,v)=\gamma_-f(t,x,-v)$  on $\Sigma_-$, which is the classical (local) inverse reflection boundary condition. Similarly, if we take $v'=v-2(v\cdot \mathbf{n}(x))\mathbf{n}(x)$, then we arrive at the classical (local) specular reflection boundary condition, see \cite{BCS97,Ca98}.
 We refer to \cite{BCS97} for possible minimal assumptions on $R$ such that \eqref{bcc} is well-defined, i.e., $R$ is nonnegative and it verifies the normalization condition as well as the reciprocity principle. Detailed discussions on the boundary conditions can be found in \cite{Ce90, CLL95, CGL01}.

 Here, we are more interested in the so-called diffuse reflection according to a Maxwellian with temperature of the thermal bath, which is nonlocal. Denote by $M_i(v)$ the Maxwellians for charged particles
\be
 M_i(v)=\frac{1}{(2\pi)^\frac{d-1}{2}\kappa_i^\frac{d+1}{2}}e^{-\frac{1}{2\kappa_i}|v|^2}, \quad i=1,...,N.\label{max}
\ee
We note that $M_i$ are the zeros of the rescaled Fokker--Planck operators $ L_{FP}^i$ given in \eqref{FKO}, i.e., $ L_{FP}^i(M_i)=0$, ($i=1,...,N$).
Then we can choose a special form of $R$ in \eqref{bcc} and propose the following boundary conditions for the distribution functions (cf. \cite{Mi10, CGL01}), which are special cases of the so-called Maxwell boundary condition \cite{M, Mi10}: for given $x\in \partial \Omega$ and $t>0$,
\bea
&& \gamma_-f_i^\varepsilon= \frac{M_i(v) }{\int_{ v\cdot \mathbf{n}(x)<0} |v\cdot \mathbf{n}(x)| M_i(v) dv} \int_{ v^*\cdot \mathbf{n}(x)>0} (\gamma_+f_i^\varepsilon) v^*\cdot \mathbf{n}(x)  dv^*,\ \ \text{on}\ \Sigma_-^x.\label{bdf}
\eea
Besides, for the electric potential $\phi^\varepsilon$, we simply impose the zero-outward electric field condition such that
\be
\nabla_x \phi^\varepsilon\cdot \mathbf{n}
=0,\quad \mbox{on}\ \partial \Omega. \label{bdphi}
\ee

In summary, below we will consider the rescaled VPFP system \eqref{f}--\eqref{P} on $(0, T)\times \Omega\times \mathbb{R}^d$ subject to boundary conditions \eqref{bdf}--\eqref{bdphi} and the initial data (depending on the parameter $\varepsilon$):
\be
f_{i}^\varepsilon(t, x, v)|_{t=0}=f_{i0}^\varepsilon(x, v).\label{ini}
\ee
 We remark that the boundary conditions \eqref{bdf} allow us to preserve mass conservation and obtain proper energy and entropy balance laws of the VPFP system \eqref{f}--\eqref{P}.  Denote by
\be
n_i^\varepsilon(t,x)=\int_{\mathbb{R}^d} f_i^\varepsilon(t,x,v) dv \quad \text{and}\quad J^\varepsilon_i=\frac{1}{\varepsilon}\int_{\mathbb{R}^d} v f_i^\varepsilon dv, \label{cur}
\ee
 the densities  as well as the current densities associated to the distribution functions for the $i$-th species, respectively.
   Multiplying \eqref{bdf} by $v\cdot\mathbf{n}(x)$ and integrating over $\Sigma_-^x$, we easily deduce the (macroscopic) boundary conditions for the fluxes such that
\be
J^\varepsilon_i\cdot \mathbf{n}=0, \quad \mbox{on}\ \partial \Omega,\label{bdJ}
\ee
which imply that all the particles that reach the boundary are reflected (no particle goes out nor enters in the domain $\Omega$) and thus the mass $\int_\Omega n_i^\varepsilon(t,x) dx$ is conserved for all time. On the other hand, in order to uniquely determine the solution  $\phi^\varepsilon$ to the Poisson equation \eqref{P} with homogeneous Neumann boundary condition \eqref{bdphi}, we require the global neutrality condition
\be
\sum_{i=1}^N z_i\int_\Omega \int_{\mathbb{R}^d} f_i^\varepsilon dv dx+ \int_\Omega D(x)dx=0.\label{gloneu}
\ee
and the zero-mean constraint $\int_\Omega \phi^\varepsilon dx=0$.

\subsection{Main result}

We first introduce the definition of renormalized solutions in the spirit of \cite{GhMa10, MT07}:

\begin{definition}\label{defre}
The set $(f_i^\varepsilon, \phi^\varepsilon)\in L^\infty(0,T; (L^1(\Omega\times \mathbb{R}^d))^N\times H^1(\Omega))$ is a renormalized solution to the VPFP system  \eqref{f}--\eqref{P} with initial and boundary conditions \eqref{bdf}--\eqref{ini}, if

(1) For all functions $\beta_i \in C^2(\mathbb{R})$, $i=1,...,N$ satisfying
 $$|\beta_i(s)|\leq C(s^\frac12+1),\quad |\beta'_i(s)|\leq C(1+s)^{-\frac12},\quad  |\beta''_i(s)|\leq C(1+s)^{-1},\quad \forall\, s\geq 0,$$
  the set $(\beta_i(f_i^\varepsilon), \phi^\varepsilon)$ is a  weak solution to the system
\bea
 &&\varepsilon\partial_t \beta_i(f_i^\varepsilon)
 + v\cdot\nabla_x \beta_i(f_i^\varepsilon)
 -\kappa_i z_i\nabla_x\phi^\varepsilon\cdot\nabla_v \beta_i(f_i^\varepsilon)
 =\frac{\zeta_i}{\varepsilon}L_{FP}^i(f_i^\varepsilon) \beta_i'(f_i^\varepsilon),\label{fb}\\
 && -\varpi \Delta_x \phi^\varepsilon=\sum_{i=1}^N z_i \int_{\mathbb{R}^d} f_i^\varepsilon dv+D(x),\label{Pb}
 \eea
 with initial data
\be
\beta_i(f_0^\varepsilon)|_{t=0}=\beta_i(f^\varepsilon_{i0})\label{inib}
\ee
and boundary conditions
\bea
&& \gamma_-\beta_i(f_i^\varepsilon)
=\frac{M_i(v)}{\int_{ v\cdot \mathbf{n}(x)<0} |v\cdot \mathbf{n}(x)| M_i(v)dv} \int_{ v^*\cdot \mathbf{n}(x)>0}  \gamma_+\beta_i(f_i^\varepsilon) v^*\cdot \mathbf{n}(x)dv^*
,\label{bdfb}\\
&& \nabla_x \phi^\varepsilon \cdot\mathbf{n}=0.
\eea
 (2) For any $\lambda>0$, the functions $\theta^i_{\varepsilon, \lambda}=(f_i^\varepsilon+\lambda \widetilde{M}_i)^\frac12$  satisfy
 \bea
 \varepsilon\partial_t \theta^i_{\varepsilon, \lambda} +v\cdot\nabla_x \theta^i_{\varepsilon, \lambda}
 -\kappa_i z_i \nabla_v \cdot(\nabla_x\phi^\varepsilon \theta^i_{\varepsilon, \lambda}) =\frac{\zeta_i}{2\varepsilon\theta^i_{\varepsilon, \lambda}} L_{FP}^i(f_i^\varepsilon)
 + \frac{ z_i\lambda \widetilde{M}_i}{2\theta^i_{\varepsilon, \lambda}}v\cdot\nabla _x\phi^\varepsilon,
  \label{fc}
   \eea
 where $\widetilde{M}_i$ are the normalized Maxwellians (comparing with \eqref{max})
 \be
 \widetilde{M}_i(v)=  \left(\frac{\kappa_i}{2\pi}\right)^\frac12 M_i(v) \quad \text{such that} \ \int_{\mathbb{R}^d} \widetilde{M}_i(v) dv=1, \quad i=1,...,N.
 \label{MM}
 \ee
 \end{definition}

 \begin{remark}
 Due to the regularity of renormalized functions $\beta_i$, the corresponding boundary conditions \eqref{bdfb} for the renormalized distribution functions make sense. We refer to \cite{Ces85, BP87, U86} (see also \cite{BCS97, Mi10}) for more detailed discussions about the traces of distribution functions on the boundary.
 \end{remark}

Next, we consider the rescaled version of the PNP system \eqref{pnpp}:
\bea
&& \partial_t n_i+\nabla_x \cdot J_i=0,\label{pnp1} \\
&&-\varpi \Delta_x \phi=\sum_{i=1}^N z_i n_i+D(x),\label{pnp3}
\eea
with density currents given by
\be
J_i=-\frac{1}{\zeta_i}\nabla_x n_i-\frac{ z_i}{\zeta_i} n_i\nabla_x \phi\label{curr}
\ee
and
subject to the following boundary conditions and initial conditions:
\bea
&& J_i\cdot\mathbf{n}=\nabla_x \phi\cdot \mathbf{n}=0, \quad \text{on}\  (0,T)\times \partial\Omega,\label{pnp4}\\
&&  n_i|_{t=0}=n_{i0}, \quad \text{in} \ \Omega.\label{pnp5}
\eea
Moreover, we require that
\be
\int_\Omega \phi dx=0 \quad \text{and} \quad \int_\Omega \left(\sum_{i=1}^N z_i n_i+D(x)\right) dx=0.\non
\ee
Then we introduce the weak formulation of the PNP system \eqref{pnp1}--\eqref{pnp5}.
\bd
\label{defpnp} We say that the set  $(n_i, \phi)$ is a weak solution to the initial boundary value problem of the PNP system \eqref{pnp1}--\eqref{pnp5}, if
\bea
n_i &\in& L^\infty(0,T; L\log L(\Omega)),\quad  \sqrt{n_i} \in L^2(0,T; H^1(\Omega)), \non\\
 \partial_t n_i &\in & L^1(0,T; W^{-1,1}(\Omega)), \qquad \phi \in L^2(0,T; H^1(\Omega)),\non
\eea
where the function space $L \log L(\Omega)$ is given by
 $$L\log L(\Omega):=\left\{n: n\geq 0, \int_\Omega n(1+|\log n|) dx<+\infty\right\}$$
and the PNP system \eqref{pnp1}--\eqref{pnp3} is satisfied in the weak sense: for any $u\in C^\infty([0,T]; C^\infty(\overline{\Omega}))$, $\psi \in L^2 (0, T; (H^1(\Omega))')$,
\bea
&& \int_\Omega n_i(t, \cdot) u(t, \cdot) dx- \int_\Omega n_{i0} u(0, \cdot)dx\non\\
&& = \int_0^t\int_\Omega n_i \partial_t u dxd\tau - \frac{1}{\zeta_i} \int_0^t \int_\Omega (\nabla_x n_i+ z_i n_i\nabla_x \phi)\cdot \nabla_x u dxd\tau,\quad t\in [0,T],
\non\\
&& \varpi\int_0^T\int_\Omega \nabla_x \phi\cdot \nabla_x \psi dxdt
=\int_0^T\int_\Omega \left(\sum_{i=1}^N z_i n_i+D(x)\right) \psi dxdt.\non
\eea
Moreover, the weak solution $(n_i, \phi)$ satisfies the following energy inequality
  \bea
&& e(t)+\sum_{i=1}^N \int_0^t \int_\Omega \frac{1}{\zeta_i} n_i\Big|\nabla \Big(\ln n_i+z_i \phi\Big)\Big|^2dxdt\leq e(0), \quad t\in [0,T],\non\\
 &&\qquad \text{with}\ \  e(t):= \int_\Omega \left(\sum_{i=1}^N n_i\ln n_i +\frac{\varpi}{2}|\nabla \phi|^2\right) dx.\non
 \eea
\ed
Now we are in a position to state the main result of this paper.

\bt \label{limit}
Let the background charge $D$ be independent of time and satisfy
  $D(x)\in L^\infty(\Omega)$. We assume that the initial data $(f_{i0}^\varepsilon, \phi_0^\varepsilon)$ satisfy the following assumptions
\bea  && f^\varepsilon_{i0}\geq 0, \quad \int_\Omega\int_{\mathbb{R}^d} f^\varepsilon_{i0} (1+ |v|^2
+|\log f^\varepsilon_{i0} |) dvdx\leq C_0,\label{nonn}\\
&& \|\phi^\varepsilon_0\|_{H^1(\Omega)}\leq C_0,\quad  \int_\Omega \phi_0^\varepsilon dx=0,
\eea
for some constant $C_0 >0$ independent of the parameter $\varepsilon$, and the global neutrality condition holds
\be
\sum_{i=1}^N z_i\int_\Omega \int_{\mathbb{R}^d} f_{i0}^\varepsilon dv dx+ \int_\Omega D(x)dx=0, \quad \forall\, \varepsilon>0.\label{gloneu0}
\ee
Let $(f_i^\varepsilon,  \phi^\varepsilon)$ be a free energy (renormalized) solution of the VPFP system
\eqref{f}--\eqref{P} with corresponding initial and boundary conditions \eqref{bdf}--\eqref{ini} (cf. Definition \ref{defre}). Then, as $\varepsilon$ tends to zero, up to a subsequence if necessary, we have the strong convergence results
\bea
&& f_i^\varepsilon(t,x,v) \to n_i(t,x)\widetilde{M}_i(v)\quad \, \text{in}\quad  L^1(0,T; L^1(\Omega\times\mathbb{R}^d)),\\
&& \phi^\varepsilon(t,x)\to \phi(t,x) \qquad \qquad\quad  \text{in}\quad L^2(0, T; W^{1,p}(\Omega)), \quad 1\leq p<2.\label{conphi1}
\eea
Moreover, $n_i^\varepsilon$ strongly converge in $L^1(0,T; L^1(\Omega))$ towards $n_i$ and $(n_i,\phi)$ is a weak solution to the PNP system \eqref{pnp1}--\eqref{pnp5} (cf. Definition \ref{defpnp}) with initial data  $n_i|_{t=0}=n_{i0}=\int_{\mathbb{R}^d} f_{i0} dv$,  such that $f_{i0}$ are the weak limits of $f^\varepsilon_{i0}$.
\et
\begin{remark}
We would like to mention that the PNP system \eqref{pnp1}--\eqref{pnp5} can also be derived from diffusion limits of other types of kinetic equations, e.g.,
the Boltzmann--Poisson system. We refer to \cite{MT07} for the one species case and we believe that their argument can also be extended to the multi-species case.
\end{remark}
 \begin{remark}\label{high}
 We remark that different types of scalings can be chosen for the VPFP system. For instance, if we assume that the drift and thermal velocities are comparable, but both are small comparing with the relaxation velocity, e.g.,
 $\nu=\mathcal{O}(1)$ and $\varepsilon<<1$, then we arrive at a different rescaled VPFP system
  \be
  \begin{cases}
   \partial_t f_i^\varepsilon + v\cdot \nabla_{x} f_i^\varepsilon - \displaystyle{\frac{\kappa_i z_i}{\varepsilon}} \nabla_{x}\phi^\varepsilon \cdot\nabla_{v} f_i^\varepsilon  = \frac{\zeta_i}{\varepsilon} L^i_{FP}(f_i^\varepsilon ), \non\\
   -\varpi  \Delta_x \phi^\varepsilon =\displaystyle{\sum_{i=1}^n z_i \int_{\mathbb{R}^d} f_i^\varepsilon  dv} +D(x).\non
   \end{cases}
 \ee
 This is usually called drift-collision balance scaling or high field scaling in the literature. Taking the hydrodynamic limit as $\varepsilon\to 0$ (the high field limit or the hyperbolic limit), the above VPFP system will lead to a first-order hyperbolic system for the density of particles coupled with the Poisson equation, cf. e.g., \cite{ACGS01,NPS01,GaS,GNPS}.
 \end{remark}

\section{Proof of Theorem \ref{limit}}
\setcounter{equation}{0}
\subsection{Uniform estimates and existence}

The  free energy of the VPFP system \eqref{f}--\eqref{P} is defined as follows
\bea
\mathcal{E}(t)=\sum_{i=1}^N \int_\Omega\int_{\mathbb{R}^d} \left(\frac{|v|^2}{2\kappa_i} f_i^\varepsilon+\mathcal{H}( f_i^\varepsilon)\right)dvdx + \frac{\varpi}{2}\int_\Omega |\nabla_x \phi^\varepsilon|^2dx,
\label{freeE}
\eea
where the function $\mathcal{H}$ takes the form $\mathcal{H}(s)=s\log s$ for $s\geq 0$.
The entropy productions of the VPFP system are given by
\bea
  \mathcal{D}^i(w)&=& \int_\Omega \int_{\mathbb{R}^d} ( v\sqrt{w}+2\kappa_i\nabla_v \sqrt{w})^2 dvdx\non\\
  &=& 4 \int_\Omega \int_{\mathbb{R}^d}
  \left|\nabla_v \sqrt{w e^{\frac{1}{2\kappa_i}|v|^2}}\right|^2 e^{-\frac{1}{2\kappa_i}|v|^2}dvdx, \quad i=1,...,N.\label{enpr}
\eea
Moreover, we introduce the Darroz\`es--Guiraud information on the boundary (cf. e.g., \cite{DG66}) such that
\bea
&&\mathcal{I}^i(w)=\int_{\Sigma^x_+} \mathcal{H}\left(w\right) d\mu^i_x- \mathcal{H} \left( \int_{\Sigma^x_+} w d\mu^i_x\right),\quad i=1,...,N,\non
\eea
 where $d\mu^i_x(v)=M_i(v)|v\cdot \mathbf{n}(x)|dv$  are  probability measures on $\Sigma^x_\pm$ by the particular choice of the normalized Maxwellians $M_i$ (cf. \eqref{max}).

First, we derive the energy dissipation property of the VPFP system \eqref{f}--\eqref{P} with initial and boundary conditions \eqref{bdf}--\eqref{ini}.

 \bp[Energy dissipation]  \label{energy}
 The renormalized solution of the VPFP system \eqref{f}--\eqref{P} with described initial data and boundary conditions satisfies
 \be
 \partial_t n_i^\varepsilon+\nabla_x \cdot J_i^\varepsilon=0, \label{aa1}
\ee
where $n_i^\varepsilon$ and $J_i^\varepsilon$ are given in \eqref{cur}. Moreover, the following dissipative energy inequality holds
\bea
&& \mathcal{E}(t)
     + \frac{1}{\varepsilon^2} \sum_{i=1}^N \frac{\zeta_i}{\kappa_i } \int_0^t\mathcal{D}^i(f_i^\varepsilon) ds
     + \frac{1}{\varepsilon} \sum_{i=1}^N   \int_0^t \int_{\partial\Omega} \mathcal{I}^i\left(\frac{\gamma_+ f^\varepsilon_{i}}{M_i(v)}\right) dS ds \non\\
&\leq & \mathcal{E}(0), \quad \forall\, t\geq 0. \label{freein}
\eea
 \ep
 \begin{proof}
 We just present a formal calculation which leads to \eqref{freein}. For $i=1,...,N$, multiplying the $i$-th equation in \eqref{f} of the VPFP system by $\frac{1}{2}|v|^2$ and integrating the result with respect to $v$ and $x$, we get
\bea
&&\frac{d}{dt} \int_\Omega \int_{\mathbb{R}^d} \frac{1}{2}|v|^2 f_i^\varepsilon dvdx
  + \int_\Omega \int_{\mathbb{R}^d}  \frac{1}{2\varepsilon}|v|^2 v\cdot\nabla_x f_i^\varepsilon dvdx\non\\
&& \ \ -\int_\Omega \int_{\mathbb{R}^d} \frac{\kappa_i z_i}{2\varepsilon} |v|^2\nabla_x\phi^\varepsilon\cdot\nabla_v f_i^\varepsilon dvdx\non\\
&=& \int_\Omega \int_{\mathbb{R}^d} \frac{\zeta_i}{2\varepsilon^2}|v|^2 L_{FP}^i( f_i^\varepsilon) dvdx,\non
\eea
integrating by parts, we see that
\be
 \int_\Omega \int_{\mathbb{R}^d}  \frac{1}{2\varepsilon}|v|^2 v\cdot\nabla_x f_i^\varepsilon dvdx
 =\frac{1}{2\varepsilon} \int_{\partial\Omega} \int_{\mathbb{R}^d} (v\cdot\mathbf{n}) |v|^2 \gamma f_i^\varepsilon dvdS,\non
 \ee
 \bea
   -\int_\Omega \int_{\mathbb{R}^d} \frac{1}{2\varepsilon} |v|^2\nabla_x\phi^\varepsilon\cdot\nabla_v f_i^\varepsilon dvdx
   &=&
 \frac{ 1}{\varepsilon} \int_\Omega \int_{\mathbb{R}^d} (v\cdot\nabla_x\phi^\varepsilon) f_i^\varepsilon dvdx\non\\
 &=& - \int_\Omega \phi^\varepsilon\nabla_x \cdot J^\varepsilon_i dx + \int_{\partial \Omega} \gamma\phi^\varepsilon J^\varepsilon_i \cdot \mathbf{n} dS\non\\
 &=& \int_\Omega \phi^\varepsilon\partial_t n_i^\varepsilon dx,\non
 \eea
 and
 \be
   \int_\Omega \int_{\mathbb{R}^d} \frac{1}{2\varepsilon^2}|v|^2 L_{FP}^i( f_i^\varepsilon) dvdx=-\frac{1}{\varepsilon^2}\int_\Omega \int_{\mathbb{R}^d}(vf_i^\varepsilon+ \kappa_i \nabla_v f_i^\varepsilon)\cdot v dvdx.\non
\ee
As a result, for $i=1,...,N$ we obtain  that
 \bea
 &&\frac{d}{dt} \int_\Omega \int_{\mathbb{R}^d} \frac{1}{2\kappa_i}|v|^2 f_i^\varepsilon dvdx + z_i\int_\Omega \phi^\varepsilon\partial_t n_i^\varepsilon dx\non\\
 &=& -\frac{1}{2\kappa_i \varepsilon} \int_{\partial\Omega} \int_{\mathbb{R}^d} (v\cdot\mathbf{n}) |v|^2 \gamma f_i^\varepsilon dvdS
  -\frac{\zeta_i}{\kappa_i\varepsilon^2}\int_\Omega \int_{\mathbb{R}^d}( vf_i^\varepsilon+\kappa_i \nabla_v f_i^\varepsilon)\cdot v dvdx.
   \label{E1}
 \eea
 Next, multiplying the $i$-th  equation \eqref{f} of the VPFP system by $\log f_i^\varepsilon$ and integrating the result with respect to $v$ and $x$, we get
 \bea
 && \frac{d}{dt}\int_\Omega \int_{\mathbb{R}^d} \mathcal{H}(f_i^\varepsilon) dvdx
 + \int_\Omega \int_{\mathbb{R}^d} \frac{1}{\varepsilon}(v\cdot\nabla_x f_i^\varepsilon ) \log f_i^\varepsilon dv dx\non\\
 &&\ \
 -\int_\Omega \int_{\mathbb{R}^d} \frac{\kappa_i z_i}{\varepsilon} (\nabla_x \phi^\varepsilon\cdot\nabla_v f_i^\varepsilon)\log f_i^\varepsilon dvdx\non\\
 &=&\int_\Omega \int_{\mathbb{R}^d} \frac{\zeta_i}{\varepsilon^2} L_{FP}^i(f_i^\varepsilon) \log f_i^\varepsilon dvdx,\non
 \eea
after integrating by parts, we see that
 \bea
  && \int_\Omega \int_{\mathbb{R}^d} \frac{1}{\varepsilon} (v\cdot\nabla_x f_i^\varepsilon ) \log f_i^\varepsilon dv dx\non\\
  &=& -\frac{1}{\varepsilon}\int_\Omega \int_{\mathbb{R}^d} (v\cdot\nabla_x f_i^\varepsilon ) dv dx+
  \frac{1}{\varepsilon}\int_{\partial\Omega}\int_{\mathbb{R}^d} (v\cdot\mathbf{n})\gamma f_i^\varepsilon\log \gamma f_i^\varepsilon dvdS\non\\
  &=&-  \int_{\partial\Omega} J_i^\varepsilon \cdot \mathbf{n} dS +
  \frac{1}{\varepsilon}\int_{\partial\Omega}\int_{\mathbb{R}^d} (v\cdot\mathbf{n})\gamma f_i^\varepsilon\log\gamma  f_i^\varepsilon dvdS\non\\
  &=& \frac{1}{\varepsilon}\int_{\partial\Omega}\int_{\mathbb{R}^d} (v\cdot\mathbf{n})\gamma f_i^\varepsilon\log\gamma  f_i^\varepsilon dvdS,\non
 \eea
 \be
 -\int_\Omega \int_{\mathbb{R}^d} \frac{1}{\varepsilon} (\nabla_x\phi^\varepsilon\cdot\nabla_v f_i^\varepsilon)\log f_i^\varepsilon dvdx=0,\non
 \ee
 and
 \be
 \int_\Omega \int_{\mathbb{R}^d} \frac{\zeta_i}{\varepsilon^2} L_{FP}^i(f_i^\varepsilon) \log f_i^\varepsilon dvdx
 =-\frac{\zeta_i}{\varepsilon^2} \int_\Omega \int_{\mathbb{R}^d}( vf_i^\varepsilon+ \kappa_i \nabla_v f_i^\varepsilon)\cdot \frac{\nabla_v f_i^\varepsilon}{f_i^\varepsilon}dvdx.\non
 \ee
 As a consequence, we get
\bea
  \frac{d}{dt}\int_\Omega \int_{\mathbb{R}^d} \mathcal{H}(f_i^\varepsilon) dvdx
 &=&-\frac{1}{\varepsilon}\int_{\partial\Omega}\int_{\mathbb{R}^d} (v\cdot\mathbf{n})\gamma f_i^\varepsilon\log\gamma  f_i^\varepsilon dvdS\non\\
 && -\frac{\zeta_i}{\varepsilon^2} \int_\Omega \int_{\mathbb{R}^d}( vf_i^\varepsilon+ \kappa_i \nabla_v f_i^\varepsilon)\cdot \frac{\nabla_v f_i^\varepsilon}{f_i^\varepsilon}dvdx.\label{E3}
 \eea
 Moreover, we see that for $f_i^\varepsilon$
 \be
   \int_\Omega \int_{\mathbb{R}^d}(vf_i^\varepsilon+ \kappa_i\nabla_v f_i^\varepsilon)\cdot \left( v+\frac{\kappa_i \nabla_v f_i^\varepsilon}{f_i^\varepsilon}\right)dvdx= \mathcal{D}^i(f_i^\varepsilon).\label{E5}
 \ee
 Due to the Poisson equation \eqref{P}, we have
 \be
  \int_\Omega \phi^\varepsilon\partial_t \left(\sum_{i=1}^N z_i n^\varepsilon_i\right) dx= -\varpi\int_\Omega \phi^\varepsilon\partial_t \Delta_x\phi^\varepsilon dx =\frac{\varpi}{2} \frac{d}{dt} \int_\Omega |\nabla_x \phi^\varepsilon|^2dx.
  \label{E6}
 \ee
Then we conclude from \eqref{E1}--\eqref{E6} that
\be
\frac{d}{dt} \mathcal{E}(t)
+\sum_{i=1}^N \frac{\zeta_i}{\kappa_i\varepsilon^2} \mathcal{D}^i(f_i^\varepsilon)+ \frac{1}{\varepsilon} \sum_{i=1}^N \int_{\partial\Omega} \int_{\mathbb{R}^d} (v\cdot\mathbf{n}) \left(\frac{1}{2\kappa_i} |v|^2+\log \gamma f_i^\varepsilon\right) \gamma f_i^\varepsilon dvdS=0.\non
\ee
Recall that $d\mu^i_x(v)=M_i(v)|v\cdot \mathbf{n}(x)|dv$ are probability measures on $\Sigma^x_\pm$ (see the definition of $M_i(v)$ \eqref{max} and \eqref{MM}). Then for the boundary terms, we can apply the  Darroz\`{e}s--Guiraud inequality \cite{DG66}, namely, thanks to \eqref{bdf}, the convexity of $\mathcal{H}(s)=s\log s$ and the Jensen inequality we deduce that (see also \cite{Mi10})
\bea
&& \int_{\partial\Omega} \int_{\mathbb{R}^d} (v\cdot\mathbf{n}) \left(\frac{1}{2\kappa_i} |v|^2+\log \gamma f_i^\varepsilon\right) \gamma f_i^\varepsilon dvdS\non\\
&=& \int_{\partial\Omega} \int_{\Sigma^x_+} \mathcal{H}\left(\frac{\gamma_+ f^\varepsilon_{i}}{M_i(v)}\right) d\mu^i_xdS
-\int_{\partial\Omega} \int_{\Sigma^x_-} \mathcal{H}\left(\frac{\gamma_- f^\varepsilon_{i}}{M_i(v)}\right) d\mu^i_xdS\non\\
&=& \int_{\partial\Omega} \int_{\Sigma^x_+} \mathcal{H}\left(\frac{\gamma_+ f^\varepsilon_{i}}{M_i(v)}\right) d\mu^i_xdS
-\int_{\partial\Omega}\mathcal{H}\left(\int_{\Sigma^x_+}\frac{\gamma_+ f^\varepsilon_{i}}{M_i(v)} d\mu^i_x\right)dS\non\\
&=& \int_{\partial\Omega} \mathcal{I}^i\left(\frac{\gamma_+ f^\varepsilon_{i}}{M_i(v)}\right) dS\non\\
&\geq& 0.\non
\eea
 As a consequence,
\be
\frac{d}{dt} \mathcal{E}(t)+ \frac{1}{\varepsilon^2} \sum_{i=1}^N \frac{\zeta_i}{\kappa_i} \mathcal{D}^i(f_i^\varepsilon)+ \frac{1}{\varepsilon} \sum_{i=1}^N \int_{\partial\Omega} \mathcal{I}^i\left(\frac{\gamma_+ f^\varepsilon_{i}}{M_i(v)}\right) dS\leq  0. \label{entro}
\ee
Integrating \eqref{entro} with respect to time, we arrive at our conclusion \eqref{freein}.
 \end{proof}
 The energy dissipation \eqref{freein} yields the following global estimates that are uniform in the parameter $\varepsilon$, which enable us to take the diffusion limit as $\varepsilon\to 0$:
\bl \label{unie}
For any $T>0$, there exists a constant $C$ depending on $C_0$, $\zeta_i$, $\kappa_i$, $\varpi$, but independent of $\varepsilon$ and $t\in [0,T]$ such that
\bea
&& \int_\Omega \int_{\mathbb{R}^d} (1+|v|^2+|\log (f_i^\varepsilon)|)f_i^\varepsilon dvdx\leq C,\non\\
&&  \int_\Omega |\nabla_x \phi^\varepsilon|^2dx \leq C,\non\\
&&  \frac{1}{\varepsilon^2}\int_0^t \mathcal{D}^i(f_i^\varepsilon) ds\leq C,\non\\
&&  \frac{1}{\varepsilon}\int_0^t \int_{\partial\Omega} \mathcal{I}^i\left(\frac{\gamma_+ f^\varepsilon_{i}}{M_i(v)}\right) dS ds\leq C.\non
\eea
The functions $f_i^\varepsilon$ are weakly relatively compact in $L^1((0,T)\times \Omega\times \mathbb{R}^d)$ and fulfill
\be
\|\nabla_v \sqrt{f_i^\varepsilon}\|_{L^2((0,T)\times\Omega\times\mathbb{R}^d)}\leq C.\non
\ee
Concerning the fluxes, we have
\bea
&& \|J^\varepsilon_i(t, \cdot)\|_{L^1(\Omega)}\leq  \frac{1}{2\varepsilon^2} \mathcal{D}^i(f_i^\varepsilon)+ \frac12\| f^\varepsilon_{i0}\|_{L^1(\Omega\times\mathbb{R}^d)}.\non
\eea
\el
\begin{proof}
The proof is similar to \cite[Propositions 5.1, 5.2, 5.3]{GhMa10}, based on the energy inequality \eqref{freein}. Since we are now dealing with the bounded domain, we do not need to estimate terms like $\int_\Omega \int_{\mathbb{R}^d} |x| f_i^\varepsilon dvdx$ as in \cite{GhMa10}. The $L^1$ weak compactness of $f_i^\varepsilon$ follows from the well-known Dunford--Pettis theorem.
\end{proof}

   We recall that the initial boundary value problem of a full Vlasov--Poisson--Fokker--Planck--Boltzmann system (subject to more general reflection boundary conditions for the distribution function but only for one species of charged particles) has been studied in the recent paper \cite{Mi10}. The author proved the existence of DiPerna--Lions renormalized solutions by using the approximation procedure in \cite{Mi00} with some  crucial trace theorems previously introduced by the same author for the Vlasov equations \cite{Mi00a} and some new results concerning weak-weak convergence (the renormalized convergence and the biting $L^1$-weak convergence). For the current case with multiple species of charged particles, the coupling between different species is somewhat weak, i.e., only via the Poisson equation. As a result, based on the energy dissipation property Proposition \ref{energy} and Lemma \ref{unie}, we are able to prove the following existence result on renormalized solutions to the VPFP system \eqref{f}--\eqref{P}, by adapting the argument in \cite{Mi10} (see also \cite{MT07, Mi00, BD95}) with minor modifications. The details are thus omitted.

  \bt[Existence of renormalized solution] \label{exe}  Suppose that the assumptions \eqref{nonn}--\eqref{gloneu0} on the initial data are satisfied.
 For arbitrary but fixed $\varepsilon>0$, the initial boundary value problem of the VPFP system \eqref{f}--\eqref{P} admits at least one (renormalized) solution $(f_i^\varepsilon, \phi^\varepsilon)$ in the sense of Definition \ref{defre}, which satisfies Proposition \ref{energy}.
 \et


\subsection{Low field limit as $\varepsilon\to 0$}

The proof of Theorem \ref{limit} mainly follows the arguments in \cite{GhMa10} for the VPFP system that concerns only one single species of particles in the whole space.
However, for the present problem involving multiple species of charged particles, we need to modify the previous argument to deal with nonlocal interactions between particles as well as the boundary conditions. In what follows, we state the essential steps and point out the possible differences in the proof.

\textit{Step 1. Strong convergence of the electric potential $\phi^\varepsilon$}.

Based on the uniform estimates in Lemma \ref{unie}, it is straightforward to argue as \cite[Propisition 3.3]{MT07} to conclude that

\bl\label{bdnphi}
The renormalized solution $(f_i^\varepsilon, \phi^\varepsilon)$ satisfies the following properties:
\begin{itemize}
\item[(1)] for $i=1,...,N$, $n_i^\varepsilon(t,x)=\int_{\mathbb{R}^d} f_i^\varepsilon(t,x,v) dv$ are weakly relatively compact in $L^1((0, T)\times \Omega)$,
\item[(2)] $ \phi^\varepsilon(t,x)$ is relatively compact in $L^2(0, T; W^{1,p}(\Omega))$ with $1\leq p<2$.
\end{itemize}
\el

Therefore, the strong convergence of $\phi^\varepsilon$ \eqref{conphi1} (up to a subsequence) is a direct consequence of Lemma \ref{bdnphi}.

\textit{Step 2. Strong convergence of the charge densities $n_i^\varepsilon$}.

Lemma \ref{bdnphi} also implies the weak compactness of densities $n_i^\varepsilon$. Indeed, we can show the convergence of density functions in the strong sense. By using the definition of renormalized solutions (cf. Definition \ref{defre}) and a velocity averaging lemma (cf. \cite[Lemma 4.2]{MT07}, also \cite{DLM91}), we are able to obtain the compactness of the densities (cf. \cite[Proposition 6.1]{GhMa10}) such that
the densities $n_i^\varepsilon$ are relatively compact in $L^1((0,T)\times \Omega)$, namely,
there exist $n_i\in L^1((0,T)\times \Omega)$ and up to a subsequence if necessary,
\bea
 && n_i^\varepsilon\to n_i,\quad  \text{in} \ L^1((0,T)\times\Omega) \ \text{and a.e. as}\ \varepsilon\to 0. \label{cden}
 \eea
 The above result and the simple inequality $(\sqrt{a}-\sqrt{b})^2\leq |a-b|$ imply that
\bea
 &&\sqrt{ n_i^\varepsilon}\to \sqrt{n_i},\quad \text{in} \ L^2((0,T)\times\Omega)\ \text{and a.e. as}\ \varepsilon\to 0.\label{cdena}
 \eea

 \textit{Step 3. Strong convergence of the distribution functions $f_i^\varepsilon$}.

 We recall the classical Csiszar--Kullback inequality (cf. \cite[Theorem 3.1, Section 4, pp. 314]{C}, see also \cite{K}) that for
all non-negative $u\in L^1(\mathbb{R}^d,d\mu)$ (where $d\mu$ is a probability measure) with $\int_{\mathbb{R}^d} u d\mu=1$, it holds
$$\|u-1\|_{L^1(\mathbb{R}^d, d\mu)}\leq 2\left( \int_{\mathbb{R}^d} (u\log u-u+1)d\mu\right)^\frac12.$$
Choose in the above inequality
$$u=\frac{f^\varepsilon_i}{n_i^\varepsilon \widetilde{M}_i(v)}, \quad d\mu=\widetilde{M}_i(v)dv,$$
which easily implies
\be
 \left(\int_\Omega\int_{\mathbb{R}^d} |f_i^\varepsilon-n_i^\varepsilon \widetilde{M}_i(v)|dvdx\right)^2
\leq 4\left(\int_\Omega n_i^\varepsilon dx\right) \int_\Omega\int_{\mathbb{R}^d} f_i^\varepsilon\log\left(\frac{f_i^\varepsilon}{n_i^\varepsilon \widetilde{M}_i(v)}\right) dvdx. \label{CK}
 \ee
 Next, we proceed to estimate the second factor in the righthand side of \eqref{CK}.  Recalling the logarithmic Sobolev inequality (cf. e.g., \cite[Corollary 4.2]{Gross})
 \bea
 && \int_{\mathbb{R}^d} |h(v')|^2\log |h(v')| d\mu(v')\non\\
 &\leq& \int_{\mathbb{R}^d} |\nabla_{v'} h(v')|^2 d\mu(v')+ \|h(v')\|_{L^2(\mathbb{R}^d, d\mu(v'))}^2\log \|h(v')\|_{L^2(\mathbb{R}^d, d\mu(v'))},
\non
\eea
where $d\mu(v')$ is the Gauss measure
$d\mu(v')=(2\pi)^{-\frac{d}{2}}e^{-\frac{|v'|^2}{2}}dv$. Making the simple change of variable $v'\to \frac{v}{\sqrt{\kappa}}$ and denoting $h_{\kappa}(v)=h(v')$, we have
$$d \mu(v') = \left(\frac{\kappa}{2\pi}\right)^\frac{d}{2} e^{-\frac{1}{2\kappa}|v|^2}dv:= d\mu_\kappa(v),\quad  \|h(v')\|_{L^2(\mathbb{R}^d, d\mu(v'))} = \|h_\kappa(v)\|_{L^2(\mathbb{R}^d, d\mu_\kappa(v))},$$
which yields that
\bea
 && \int_{\mathbb{R}^d} |h_\kappa(v)|^2\log |h_\kappa(v)| d\mu_\kappa(v)\non\\
 &\leq&
  \kappa \int_{\mathbb{R}^d} |\nabla_{v} h_\kappa(v)|^2 d\mu_\kappa(v)+ \|h_\kappa(v)\|_{L^2(\mathbb{R}^d, d\mu_\kappa(v))}^2\log \|h_\kappa(v)\|_{L^2(\mathbb{R}^d, d\mu_\kappa(v))}.\non
\eea
In the above inequality, we set
\be
\kappa=\kappa_i, \quad h_\kappa(v)=\sqrt{\frac{f_i^\varepsilon}{ \widetilde{M}_i(v)}},\quad d\mu_\kappa(v)=\widetilde{M}_i(v) dv.  \non
\ee
Then we infer from the definition \eqref{cur}  that
$ \|h_\kappa\|_{L^2(\mathbb{R}^d, d\mu_\kappa(v))}^2=n_i^\varepsilon$, which yields
\bea
&&\int_\Omega \int_{\mathbb{R}^d} f_i^\varepsilon\log\left(\frac{f_i^\varepsilon}{n_i^\varepsilon \widetilde{M}_i(v)}\right) dv dx\non\\
&=& \int_\Omega \left(2\int_{\mathbb{R}^d} |h_\kappa(v)|^2\log |h_\kappa(v)| d\mu_\kappa(v)-2\|h_\kappa(v)\|_{L^2(\mathbb{R}^d, d\mu_\kappa(v))}^2\log \|h_\kappa(v)\|_{L^2(\mathbb{R}^d, d\mu_\kappa(v))}\right)dx\non\\
&\leq& 2\kappa_i \int_\Omega \int_{\mathbb{R}^d} \left|\nabla_v \sqrt{\frac{f_i^\varepsilon}{ \widetilde{M}_i(v)}}\right|^2  \widetilde{M}_i(v) dvdx\non\\
&=& \frac{\kappa_i}{2} \mathcal{D}^i(f_i^\varepsilon).\label{logso}
\eea

As a consequence, we infer
from the entropy dissipation in \eqref{freein}, the uniform estimates in Lemma \ref{unie} and the estimates \eqref{CK} and \eqref{logso}  that
 when  $ \varepsilon\to 0$,
 \bea
 && f_i^\varepsilon- n_i^\varepsilon\widetilde{M}_i\to 0,  \quad  \text{in} \ L^1((0,T)\times\Omega\times\mathbb{R}^d) \ \text{and a.e.} \non
 \eea
 Combing the above results with the convergence result of $n_i^\varepsilon$ \eqref{cden}, we conclude that as  $ \varepsilon\to 0$
 \bea
 && f_i^\varepsilon\to n_i\widetilde{M}_i,  \quad  \text{in} \ L^1((0,T)\times\Omega\times\mathbb{R}^d) \ \text{and a.e.} \label{confg1}
 \eea
Here and below, the convergence results are always understood to be up to a subsequence.

 \textit{Step 4. Weak convergence of the fluxes $J_i^\varepsilon$}.

 We introduce the auxiliary functions
\be
r^\varepsilon_i=\frac{\sqrt{f_i^\varepsilon}-\sqrt{n_i^\varepsilon \widetilde{M}_i(v)}}{\varepsilon \sqrt{\widetilde{M_i}(v)}}, \quad i=1,...,N.\label{rrr}
\ee
In analogy to \cite[Proposition 3.4]{MT07} and \cite[Proposition 5.5]{GhMa10}, we have
\bl\label{r}
For arbitrary $T>0$, the following uniform estimates hold
 \be
 \int_0^T\int_\Omega \int_{\mathbb{R}^d}\left( |r_i^\varepsilon|^2\widetilde{M}_i+ \varepsilon|r^\varepsilon_i|^2|v|^2\widetilde{M}_i+ \sqrt{\varepsilon}|r^\varepsilon_i|^2|v|\widetilde{M}_i \right)dvdxdt\leq C,\non
 \ee
 where $C$ is a constant that may depend on $C_0$, $\zeta_i$, $\kappa_i$, $\varpi$, but independent of $\varepsilon$ and $t\in [0,T]$.
\el
 Using the expressions of $r_i^\varepsilon$  (cf. \eqref{rrr}),
 we have
 \bea
 && f_i^\varepsilon= n_i^\varepsilon \widetilde{M}_i+2\varepsilon\widetilde{M}_i\sqrt{n_i^\varepsilon}r_i^\varepsilon +\varepsilon^2|r_i^\varepsilon|^2 \widetilde{M}_i.\label{rfe}
 \eea
 Due to the simple facts $\int_{\mathbb{R}^d} v \widetilde{M}_i(v) dv =0$,  it follows from \eqref{cur}, \eqref{cdena}  and Lemma \ref{r} that as  $ \varepsilon\to 0$
 \bea
 J_i^\varepsilon&=&
 2\sqrt{n_i^\varepsilon}\int_{\mathbb{R}^d} r^\varepsilon_i v \widetilde{M}_i dv +\sqrt{\varepsilon}\int_{\mathbb{R}^d} \sqrt{\varepsilon} |r^\varepsilon_i|^2 v\widetilde{M}_i dv\non\\
 &\to& 2 \sqrt{n_i} \int_{\mathbb{R}^d} r_i v \widetilde{M}_i dv, \quad \text{weakly in}\ L^1((0,T)\times \Omega),\non
 \eea
 where $r_i$ are the weak limits of $r^\varepsilon_i$, for $i=1,...,N$.

   It remains to identify the limit function of $J^\varepsilon_i$, which can be done by  using a similar argument as in \cite[Proposition 7.2]{GhMa10}. The strong convergence of $f_i^\varepsilon$ (see \eqref{confg1})  implies that for any fixed $\lambda>0$,
   $$\theta^i_{\varepsilon, \lambda} \to \sqrt{ (n_i+\lambda)\widetilde{M}_i}, \quad \text{as} \ \varepsilon \to 0.$$
   On the other hand, it follows from \eqref{rfe}  that for any $\lambda>0$, when $\varepsilon \to 0$, we have
   \be
   \frac{\zeta_i}{2\varepsilon} L_{FP}^i(f_i^\varepsilon) =  \zeta_i L_{FP}^i\left( \widetilde{M}_i\sqrt{n_i^\varepsilon}r_i^\varepsilon + \frac{\varepsilon}{2} |r_i^\varepsilon|^2 \widetilde{M}_i\right)\to \zeta_i\sqrt{n_i} L_{FP}^i(r_i\widetilde{M}_i).\non
   \ee
   As a consequence, in the renormalized formula \eqref{fc}, first for any fixed $\lambda>0$ passing to the limit as $\varepsilon\to 0$ and then letting $\lambda\to 0$, we obtain that
  \be
   \left(\nabla_x \sqrt{n_i}+\frac{  z_i}{2}\nabla_x \phi\sqrt{n_i}\right)\cdot v\widetilde{M}_i=\zeta_i L_{FP}^i(r_i\widetilde{M}_i),\label{ffr}
  \ee
  where $\phi$ is the limit of $\phi^\varepsilon$ (recall \eqref{conphi1}).

   On the other hand, it follows from \cite[Proposition 3.1]{GhMa10} that $\chi_j=-v_j\widetilde{M}_i$ ($i=1,...,N$, $j=1,...,d$) is the unique solution to the equation $L_{FP}^i \chi_j=v_j\widetilde{M}_i$ in $R(L_{FP}^i)\cap D(L_{FP}^i)$, where
  \bea
   L^2_{\widetilde{M}_i}(\mathbb{R}^d)&=&L^2(\mathbb{R}^d; \widetilde{M}_i^{-1}dv),\non\\
   R(L_{FP}^i)&=& \left\{f\in L^2_{\widetilde{M}_i}(\mathbb{R}^d): \int_{\mathbb{R}^d} f(v) dv=0\right\},
  \non\\
   D(L_{FP}^i)&=& \left\{f\in L^2_{\widetilde{M}_i}(\mathbb{R}^d): \nabla_v\cdot\Big(e^{-\frac{1}{2\kappa_i}|v|^2}\nabla_v( e^{\frac{1}{2\kappa_i}|v|^2} f)\Big)\in L^2_{\widetilde{M}_i}(\mathbb{R}^d)\right\}.\non
  \eea
  Since $-L_{FP}^i$ is a self-adjoint operator on $ L^2_{\widetilde{M}_i}(\mathbb{R}^d)$, using \eqref{ffr}, we have
  \bea
  J_i&=& 2\sqrt{n_i}\int_{\mathbb{R}^d} r_i v\widetilde{M}_i dv\non\\
     &=& 2\sqrt{n_i}\int_{\mathbb{R}^d} (r_i \widetilde{M}_i) L_{FP}^i(- v\widetilde{M}_i) \widetilde{M}_i^{-1} dv\non\\
     &=&  2 \sqrt{n_i}\int_{\mathbb{R}^d} L_{FP}^i(r_i \widetilde{M}_i) (-v  \widetilde{M}_i) \widetilde{M}_i^{-1} dv\non\\
     &=& \frac{2}{\zeta_i}\sqrt{n_i}\int_{\mathbb{R}^d}\left[(\nabla_x \sqrt{n_i}+\frac{ z_i}{2}\nabla_x \phi\sqrt{n_i})\cdot v\widetilde{M}_i\right](-v  \widetilde{M}_i) \widetilde{M}_i^{-1} dv\non\\
     &=& -\frac{2}{\zeta_i} \sqrt{n_i} \left(\int_{\mathbb{R}^d} v\otimes v \widetilde{M}_i dv\right) \left(\nabla_x \sqrt{n_i}+\frac{z_i}{2}\nabla_x \phi\sqrt{n_i}\right)\non\\
     &=& -\frac{2}{\zeta_i} \sqrt{n_i}\left(\nabla_x \sqrt{n_i}+\frac{ z_i}{2}\nabla_x \phi\sqrt{n_i}\right).\non
  \eea
  where we use the fact that $\int_{\mathbb{R}^d} v\otimes v \widetilde{M}_i dv=\mathbb{I}$.
  Therefore, we can see that as $\varepsilon\to 0$
 \be
  J_i^\varepsilon \to J_i:=-\frac{2}{\zeta_i}\sqrt{n_i} \left(\nabla_x\sqrt{n_i}+\frac{ z_i}{2}\nabla_x \phi \sqrt{n_i}\right).\label{Jf}
 \ee
 in the distribution sense.

 \textit{Step 5. Passage to the limit in the PDE system}.

 In order to recover the PNP system \eqref{pnp1}--\eqref{pnp5}, we state a regularity result for the density functions $n_i$ in the spirit of \cite[Lemma 7.1]{MT07}
 \bl
 Let $\Omega$ be a smooth bounded and open set in $\mathbb{R}^d$. Assume $n_i$ are positive functions belonging to $L^\infty(0, T; L^1(\Omega))$ and $\phi \in L^2(0, T; H^1(\Omega))$ that satisfy
 \bea
 &&\nabla_x \sqrt{n_i}+\frac{ z_i}{2}\nabla_x\phi \sqrt{n_i}=G_i\in L^2(0,T; L^2(\Omega)),\quad i=1,...,N,\label{Gf}\\
 &&-\varpi \Delta_x \phi=\sum_{i=1}^N z_i n_i+D(x). \non
 \eea
 Then we have
 \be
 \sqrt{n_i}\in L^2(0, T; H^1(\Omega)), \quad  \sum_{i=1}^N z_i n_i \in L^2(0, T; L^2(\Omega)),\quad  \nabla_x \phi \sqrt{n_i} \in L^2(0, T; L^2(\Omega)).\non
 \ee
 \el
 \begin{proof}
 As in \cite[Corollary 3.2]{MT07}, we take $\beta_\delta(s)=\delta^{-1}\beta(\delta s)$ where $\beta\in C^\infty(\mathbb{R})$ satisfying $\beta(s)=s$ for $-1\leq s\leq 1$, $0\leq \beta'(s)\leq 1$ for $s\in \mathbb{R}$ and $\beta(s)=2$ for $|s|\geq 3$. Then we renormalize the equations \eqref{Gf} for $\sqrt{n_i}$ such that
  \be
 \nabla_x \beta_\delta (\sqrt{n_i})+\frac{ z_i}{2}\nabla_x\phi \beta'_\delta(\sqrt{n_i})\sqrt{n_i}=G_i\beta'_\delta(\sqrt{n_i})\in L^2(0,T; L^2(\Omega)).\label{regg1}
 \ee
For any $\delta>0$, due to our choice of $\beta$ and the given regularity for $\nabla_x \phi$, we have
\be
 \|\nabla_x\phi \beta'_\delta(\sqrt{n_i})\sqrt{n_i}\|_{L^2(0,T; L^2(\Omega))}\leq \frac{3}{\delta} \|\nabla_x\phi\|_{L^2(0,T; L^2(\Omega))},\non
\ee
which together with \eqref{regg1} implies that $\nabla_x \beta_\delta (\sqrt{n_i})\in L^2(0,T; L^2(\Omega))$.
   Then we can take $L^2$ norm on both sides of the equations \eqref{regg1}, summing up with respect to $i=1,...,N$, we have
 \bea
 &&  \sum_{i=1}^N \|\nabla_x \beta_\delta (\sqrt{n_i})\|^2_{L^2(0,T; L^2(\Omega))}+\sum_{i=1}^N \frac{ z_i^2}{4}\|\nabla_x\phi \beta'_\delta(\sqrt{n_i})\sqrt{n_i}\|^2_{L^2(0,T; L^2(\Omega))}\non\\
 &&\  \ +\sum_{i=1}^N \int_0^T\int_\Omega \left[z_i\beta'_\delta(\sqrt{n_i})\sqrt{n_i} \nabla_x \beta_\delta (\sqrt{n_i})
  \right] \cdot \nabla_x\phi dx dt\non\\
 &\leq & \sum_{i=1}^N  \|G_i\|^2_{L^2(0,T; L^2(\Omega))},\non
 \eea
 where the right-hand side is independent of $\delta$. For the crossing term on the left hand side, using integration by parts, we have
 \bea
 && \sum_{i=1}^N \int_0^T\int_\Omega \left[ z_i\beta'_\delta(\sqrt{n_i})\sqrt{n_i}\nabla_x \beta_\delta (\sqrt{n_i})
   \right] \cdot \nabla_x\phi dx dt\non\\
  &=& \sum_{i=1}^N \int_0^T\int_\Omega\nabla_x \left[ z_i \tilde{\beta}_\delta(\sqrt{n_i})\right]\cdot \nabla_x \phi dxdt\non\\
  &=& \frac{1}{\varpi}\int_0^T\int_\Omega  \sum_{i=1}^N z_i \tilde{\beta}_\delta(\sqrt{n_i})\cdot \left(\sum_{i=1}^N z_in_i+ D(x)\right)  dxdt\non
 \eea
 where $\tilde{\beta}$ satisfies
 $$\tilde{\beta}(s)=\int_0^s \tau \beta'(\tau)^2 d\tau, \quad \tilde{\beta}_\delta(s)=\delta^{-2}\tilde{\beta}(\delta s), \quad \tilde{\beta}_\delta(s) \to \frac{s^2}{2}, \quad \text{as}\ \delta\to 0.$$
  Let $\delta\to 0$, we infer from the above estimates that
  \bea
 && \sum_{i=1}^N \|\nabla_x \sqrt{n_i}\|^2_{L^2(0,T; L^2(\Omega))}+\sum_{i=1}^N \frac{ z_i^2}{4}\|\nabla_x\phi\sqrt{n_i}\|^2_{L^2(0,T; L^2(\Omega))}\non\\
 &&\ \ + \frac{1}{2\varpi} \int_0^T\int_\Omega \left(\sum_{i=1}^N z_i n_i\right)^2 dxdt\non\\
 &\leq & \sum_{i=1}^N \|G_i\|^2_{L^2(0,T; L^2(\Omega))}+ \frac{1}{2\varpi}\left|\int_0^T\int_\Omega  D(x) \left(\sum_{i=1}^N z_i n_i\right) dxdt\right|\non\\
 &\leq & \sum_{i=1}^N  \|G_i\|^2_{L^2(0,T; L^2(\Omega))} +\frac{1}{4\varpi}\int_0^T\int_\Omega | D(x)|^2+  \left(\sum_{i=1}^N z_i n_i\right)^2 dxdt,\non
 \eea
 which easily yields the required regularity estimate. The lemma is proved.
 \end{proof}
 Finally, using the above regularity lemma and the convergence result \eqref{Jf}, we are able to write the currents $J_i$ as in \eqref{curr}. Then we can pass to the limit as $\varepsilon\to 0$ in the weak form of equations \eqref{aa1} as well as in the Poisson equation \eqref{Pb} to conclude that the limit functions $(n_i, \phi)$ satisfy the rescaled PNP system \eqref{pnp1}--\eqref{pnp5}.

  The proof of Theorem \ref{limit} is complete.

\section*{Acknowledgement} H. Wu was partially supported by NSF of China 11371098 and ``Zhuo Xue" program of Fudan University. C. Liu was partially supported by NSF grants DMS-1109107, DMS-1216938 and DMS-1159937.


\end{document}